\documentclass[12pt]{article}
\usepackage{amsmath,amsthm,amssymb, amstext, amsopn, amsxtra}
\begin{document}
 \begin{center}
\Large{A short proof of Grinshpon's theorem}\\[2mm]
\small{Dinesh Khurana} \\[2mm]
\footnotesize{Indian Institute of Science Education and Research, Mohali} \\ 
\footnotesize{MGSIPA Transit Campus, Sector 26 Chandigarh-160019, India}\\
{\tt dkhurana@iisermohali.ac.in}
\end{center}

In [G, Theorem], Grinshpon proved that if $S$ is a commutative subring of a ring $R$ and $A\in M_n(S)$ is invertible in $M_n(R)$, then $det(A)$ is invertible in $R$. Grinshpon's theorem
is immediate from the following result by taking $B = A^{-1}$ and using the fact that a square matrix over a commutative ring is invertible
iff its determinant is invertible. For $X\in M_n(R)$, let
$C(X) = \{Y \in M_n(R) : XY = YX\}$.\\[3mm]
{\bf Theorem.} {\em Let $A$, $B \in M_n(R)$ be such that $C(A) = C(B)$. If the entries of $A$ commute with each other, then there exists a commutative subring $T$ of $R$ such that $A, B \in M_n(T)$. }\\[3mm]
{\bf Proof.} Let $I$ denote the identity of $M_n(R)$. If $a$ is any entry of $A$, then $aI \in C(A) = C(B)$. Thus every entry of $A$ commutes with every entry of $B$. So if $b$ is any entry of $B$, then $bI \in C(A) = C(B)$ implying that the entries of $B$ commute with each other also. Thus the subring $T$ generated by all entries of $A$ and $B$ is commutative. \qed \\[3mm]
{\bf References}\\[3mm]
[G] M. Grinshpon, {\em Invertibility of matrices over subrings}, communications in Algebra 36 (2008), 2619-2624.\\[3mm]

\end{document}